
    \def\e{\epsilon}  \def\g{\gamma}      \def\om{\omega}   \def\s{\sigma}  \def\S{\Sigma}

\def\cl{{\rm {cl}}}                 \def\HD{{\rm{HD}}}  
\def\Arg{{\rm {Arg}}}           

  \def\a{\alpha}        \def\La{\Lambda}
\def\z{\zeta}
\def\e{\varepsilon}      \def\b{\beta}
\def\f{\varphi}            \def\d{\delta}

\def\C{I\!\!\!\!C}     \def\D{I\!\!D}    \def\Z{Z\!\!\!Z}
\def\R{I\!\!R}          \def\T{{\cal T}}   \def\B{{\cal B}}
\def\Crit{{\rm {Crit}}}     

\def\and{{\rm and}}   \def\ov{\overline}   
 \def\dist{{\rm {dist}}}     \def\length{{\rm {length}}}
\def\h{{\rm h}}                   \def\supp{{\rm {supp}}}
\def\Vol{{\rm {Vol}}}        \def\diam{{\rm {diam}}}
\def\Const{{\rm {Const}}}   \def\deg{{\rm {deg}}}

{\bf Density of periodic sources in the boundary of a basin of attraction for iteration of holomorphic maps,  geometric coding trees technique}

\

\centerline {by F. Przytycki\footnote{*}{Supported by Polish KBN 
Grants 210469101 "Iteracje i Fraktale" and 210909101 
"...Uklady Dynamiczne".} and A. Zdunik*}

\


\

{\bf Abstract.}  {\it We prove that  if A is the basin of 
immediate attraction to a periodic attracting or parabolic 
point for a rational map f on the Riemann sphere, then 
periodic points in the boundary of A are dense in this 
boundary. To prove this in the non simply-connected or  
parabolic  situations we prove a more abstract, geometric 
coding trees version.}

\

{\bf Introduction} 

\

Let $f:\ov\C\to \ov\C$ be a rational map of the Riemann sphere $\ov\C$. Let $J(f)$ denote its Julia set. 
We say a periodic point $p$ of period $m$ is attracting (a sink) 
if $|(f^m)'(p)|<1$,
repelling (a source) if $|(f^m)'(p)|>1$
 and parabolic if $(f^m)'(p)$ is a root of 
unity. We say that $A=A_p$ is the immediate basin of 
attraction to a sink or a parabolic point $p$ if $A$ is a component of $\ov\C\setminus 
J(f)$ such that $f^{nm}|_A\to p$ as $n\to\infty$ and 
$p\in A_p$ in the case $p$ is attracting, $p\in\partial A$
in the case $p$ is parabolic.

\

 We shall prove the following fact asked by G. Levin:

\

{\bf Theorem A.}  If $A$ is the basin of immediate attraction for a periodic 
attracting or parabolic point  for a rational map $f:\ov\C\to\ov\C$  then 
periodic points contained in $\partial A$ are dense in $\partial A$.

\

 A classical Fatou, Julia  theorem 
says that periodic sources are dense in $J(f)$. However 
these  periodic sources could only converge to $\partial A$, 
not being in $\partial A$.

\

The density of periodic points in Theorem A immediately implies the density of periodic sources because for every rational map there are only finitely many periodic points not being sources and Julia set has no isolated points.

 An idea of a proof of Theorem A using Pesin theory and Katok's 
proof of density of periodic points [K]
saying that $f^{-n}(B(x,\e))\subset B(x,\e)$ for some branches of $f^{-n}$, is also too crude. The matter is that 
the resulting fixed point for $f^n$ in $B(x,\e)$ could be outside $\partial A$. However this gives an idea for a 
correct proof.
We shall consider points in $\partial A$ together with 
"tails", some curves in $A$ along which these points are accessible. (We say $x\in \partial A$ is {\it accessible} from 
$A$ if  there exists a continuous 
curve  $\g:[0,1]\to \ov\C$ such that $\g([0,1))\subset 
A$  and $\g(1)=x$. We say then also that $x$ is accessible along $\g$.)

\

Thus proving Theorem A we shall prove in fact  something stronger:

{\bf Complement to Theorem A.} Periodic points in $\partial      A$  accessible 
from $A$ along $f$-invariant finite length curves, are dense in $\partial A$.

\

If $f$ is a polynomial (or polynomial-like) then it follows automatically that these periodic points are accessible along external rays. See [LP] for the proof 
and for the definition of external rays in the case $A$ is 
not simply-connected.

\

 It is an open problem whether all periodic sources 
in $\partial A$ are accessible from $A$, see [P3] for a discussion of this and
related problems. It was proved that this is so  
in the case $f$ is a polynomial and $A$ is the basin of 
attraction to $\infty$ in [EL], [D] and later in [Pe], [P4] in more general situations: for $f$ any rational function and 
$A$ a completely invariant (i.e.  $f^{-1}(A)=A$) basin of 
attraction to a sink or a parabolic point.

 \

The paper is organised as follows: In Section 1 we shall prove Theorem 1 directly
in the case of $A$ simply-connected, $ p$ attracting. In Section 2 we shall
introduce a more general point of view: geometric coding trees, studied and 
exploited already in [P1], [P2], [PUZ] and [PS], and 
formulate and prove 
Theorems B and C in the trees setting, which easily yield Theorem A. 

\

\

{\bf Section 1.  \  Theorem A in the case of a  simply-connected $A$ and $p$ 
attracting.}

\

Here we shall prove Theorem A assuming that $A$ is simply-connected and $p$ is attracting.

First let us state Lemma 1 which belongs to Pesin's 
Theory. 

\

{\bf Lemma 1.}  Let $(X,\cal F , \nu)$ be a measure space with a measurable automorphism $T:X\to X$  .
Let $\mu$ be an ergodic $f$-invariant measure on a 
compact set $Y$ in the Riemann sphere, for $f$ a holomorphic mapping from a neighbourhood of $Y$ 
to $\ov\C$ keeping $Y$ invariant, 
with positive Lyapunov exponent i.e. 
$\chi_\mu (f):=\int \log |f'| d\mu \ >0$.
Let $h:X\to Y$ be a measurable mapping such that 
$h_* (\nu)=\mu \ \ \and \ \ h\circ T=f\circ h$ a.e. .

Then for $\nu$-almost every $x\in X$ there exists 
$r=r(x)>0$ such that  univalent branches 
$F_n$ of $f^{-n}$ on $B(h(x),r)$ for $n=1,2,...$ for which  
$F_n(h(x))=h(T^{-n}(x))$, exist.
Moreover for an arbitrary $\exp (-\chi_\mu(f))
  <  \lambda < 1$ (not depending on $x$) and a constant $C=C(x)>0$ 
$$
|F_n'(h(x))| < C\lambda^n \ \ \  \and  \ \ \ 
{|F_n'(h(x))\  \over  |F_n'(z)|} <C
$$
for every $z\in B(h(x),r), \ n>0$, 
(distances  and derivatives in the Riemann metric on 
$\ov\C$).

Moreover $r$ and $C$ are measurable functions of $x$.

\

Let $R:\D\to A_p$ be a Riemann mapping such that 
$R(0)=p$. Define $g:=R^{-1}\circ f \circ R$ on $\D$. 
We know that $g$ extends holomorphically to a neighbourhood of $\cl A$ and is expanding on $\partial A$,
see [P2]. (In fact $g$ is a finite Blaschke product, because we assume in this section that $f$ is defined on the whole $A$, see [P1]. However we need only the assumption that $f$ is defined on a neighbourhood of $\partial A$ as in [P2].)

For every $\z\in \partial\D$ every $0<\a<\pi /2$ and 
every  $\rho >0$ consider the cone

$${\cal C}_{\a,\rho}(\z):=\{z\in\D : |\Arg\z - \Arg (\z -z)| < \a , \ 
|\z - z| <\rho \}
$$

 In the sequel we shall need the following simple 

 \

{\bf Lemma 2.}  There exist $\rho_0>0, C>0$  and 
$0<\a_0<\pi/2$ such that for every $\rho \le \rho_0, \ n\ge 0 , \ \z\in \partial\D$ and every 
branch $G_n$ of $g^{-n}$ on the disc $B(\z,\rho_0)$   the following inclusion holds:
$$
G_n(\{z\in \D: z=t\z, 1-t < \rho \}) \subset 
{\cal C}_{\a_0, C\rho}(G_n(\z))
$$

{\bf Remark. } Considering an iterate of $f$ and $g$ we 
can assume that $C=1$, because above we can write 
in fact ${\cal C}_{\a_0,C\xi^n\rho}$  for a number $0 < \xi < 1$.

\

{\bf Proof of Theorem A in the case of a simply-
connected basin of a sink.}

Keep the notation of this section: $A$ the basin of attraction to a fixed point, a sink $p$,  a Riemann mapping $R:\D\to A$ and $g$ the  pull-back of $f$ extended beyond $\partial\D$, just a finite Blaschke product. 

Consider $\mu:=\ov R_*(l)$, where $\ov R$ denotes the radial limit of $R$ and $l$ is the normalized length measure on $\partial\D$. In fact $\mu$ is the harmonic measure on $\partial A$ viewed from $p$. This measure is 
ergodic $f$-invariant and $\chi_\mu (f)=\chi_l(g)>0$, see [P1, P2].  Also supp $\mu=\partial A$. 

Indeed for every $\e>0, \ x\in \partial A$ and $x_n\in A$  
such that $x_n\to x$ we have for harmonic measures:  
$\om(x_n,B(x,\e)) \to 1\not= 0$. But the measures 
$\om(p,\cdot)$ and $\om(x_n,\cdot)$ are equivalent 
hence $\om(p, B(x,\e))>0$. 

We shall not use anymore the assumption $\mu$ is a 
harmonic measure, we shall use only  the 
abovementioned properties.
 
\

  From the existence of a nontangential limit $\ov R$ of $R$ a.e. 
[Du] it follows easily that  for an arbitrary $\e>0$ and 
$0<\a<\pi/2$ and $\rho >0$ there 
exists $K_\e \in \partial\D$ such that $l(K_\e)\ge 1-\e$ satisfying 
$$
R(z)\to \ov R(\z)\ \  \hbox{uniformly as}\ \ \  z\to \z, z\in 
{\cal C}_{\a,\rho}(\z)
$$
Namely for every $\d_1>0$ there exists $\d_2>0$  
such that for every $\z\in\partial \D$ if $z\in {\cal C}_{\a,\d_2}$ 
then  $\dist (R(z), 
\ov R(\z)) <\d_1$, distance in the Riemann metric on                                                                                                                                                                                     
$\ov\C$.                   

\

Consider the inverse limit  (natural extension in Rohlin 
terminology [Ro]) 
$(\tilde{\partial\D} ,\tilde l, \tilde g)$ of  
$(\partial\D, l, g).$ 

Denote the standard projection of $\tilde{\partial\D}$ 
to $\partial\D$  (the zero coordinate) by $\pi_0$.

Due to Lemma 1 applied to $(\tilde{\partial\D}, \tilde l, 
borel)$  the automorphism $\tilde g$ the map $h=\ov R \circ 
\pi_0$ and $Y=\partial A, \  f$ our rational map, there exist constants $C,r>0$ this time not dependent on $x$, 
and a measurable 
set  $\tilde K\subset\tilde{\partial\D} $ such that
$\tilde l(\tilde K) \ge 1-2\e, \ \  \tilde K \subset \pi^{-
1}_0(K_\e)$ and for every $g$-trajectory $(\z_n)\in 
\tilde K$ the assertion of Lemma 1 with the 
constants $C$ and $r$ holds. 

Let $t=t(r)$ be such a number that for every $\z\in K_\e$ 
and $z\in K_{\a,t}$  we have 
$$
\dist (R(z), \ov R(\z)) < r/3 \eqno (1)
$$

We additionally assume that $t<\rho_0$ from 
Lemma 2. Also $\a$ is that from Lemma 2.

\

By Poincar\'e Recurrence Theorem for $\tilde g$ 
for a.e. trajectory $(\z_n)\in \tilde K$ there exists a 
sequence $n_j \to \infty$ as $j\to\infty$ such that
$$
\z_{-n_j} = \pi_0 \tilde g^{-n_j}((\z_n)) \to \z_0. \eqno (2)
$$
and $\tilde g^{-n_j}((\z_n))\in \tilde K$ hence 
$$
\z_{-n_j} \in K_\e \eqno (3)
$$
Indeed, we can take a sequence of finite partitions 
${\cal A}_j$ of $\pi_0(\tilde K)$ such that the maximal 
diameters of sets of ${\cal A}_j$ converge to 0 as 
$j \to \infty$.  Almost every $(\z_n)\in \tilde K$ is in  
$\bigcap_j \pi^{-1}_0(A_j)$ where $A_j \in {\cal A}_j$ and 
there exists $n_j$ such that 
$\tilde g^{-n_j}((\z_n))  \in  \pi^{-1}_0(A_j)$

\

For a.e. $(\z_n) \in \tilde K$ fix $N=N((\z_n)) $ such that  
$$
\z_{-N} \in B(\z_0, t(r) \sin\a) \eqno (4)
$$
arbitrarily large.

Denote by $G_N$ the branch of $g^{-N}$ such that 
$G_N(\z_0)=\z_{-N}$. By Lemma 2 
$G_N((\tau\z_0) \in {\cal C}_{\a, t}(\z_N)$  for every 
$1-t<\tau<1$

By (4) there exists $1-t <\tau_0 < 1$ such that 
$\tau_0\z_0 \in {\cal C}_{\a,t}(\z_{-N})$, see Fig 1. :

\

\

\

\

\

\

\

\

\

\

\

\centerline {Figure 1}

\

Due to (3) we can apply (1) for $\z_{-N}$.  Thus by (1) 
applied to $z=\tau_0\z_0, \ \z=\z_0 \ \ \and \ \ 
\z=\z_{-N}$ we obtain 
$$
\dist (\ov R(\z_{-N}) , \ov R(\z_0))  <  {2\over 3}r .
$$

So, if $N$ has been taken large enough,  we obtain  
by Lemma 1 for the branch $F_N$ of $f^{-N}$ discussed 
in the statement of Lemma 1
$$
F_N(B(\ov R (\z_0), r)) \subset B(\ov R(\z_{-N}), r/3) \subset  
B(\ov R(\z_0), r), \eqno (5)
$$
see Fig. 2

\

\

\

\

\

\

\

\

\

\

\centerline {Figure 2}

\
 
Moreover $F_N$ is a contraction, i.e. 
$|(F_N|_{B(\ov R(\z_0), r)})'| < C\lambda^N < 1$.

The interval $I$ joining $\tau_0\z_0$ with $G_N((1-t)\z_0)$ is 
in ${\cal C}_{\a,t}(\z_{-N})$, hence 
$$
R(I) \subset  B(\ov R(\z_{-N}), r/3) \subset  B(\ov  R(\z_0), r)
$$
By the definitions of $F_N, G_N$ we have $
\ov R \circ G_N = F_N \circ \ov R $ at $\z_0$.  To 
prove this equality on $[(1-t)\z,\z]$ we must know 
that for $f^{-N}$ we have really the branch $F_N$.
But this is the case because the maps involved 
are continuous on the domains under consideration 
and $[(1-t)\z_0,\z_0]$ is connected. So 
$$ F_N(R(1-t)\z_0) = RG_N((1-t)\z_0) \eqno (6)
$$

Let $\g$ be the concatenation of the curves $R([(1-t)\z_0,\tau\z_0])$ and 
$R(I)$. By (6) it  joins $R((1-t)\z_0)$ with 
$F_N(R((1-t)\z_0))$ and it is entirely in  $B(\ov R(\z_0), r)$.  One  end  
$a$ of the curve $\Gamma$ 
being the concatenation of $\g, F_N(\g), F_N^2(\g),...$  is in $\partial A$ and 
is periodic  of period $N$, \  ($\Gamma $ makes sense 
due to  (5)). Moreover 
$$
\length(\Gamma) \le \sum_{n\ge 0}C\lambda^n \length 
\g <\infty
$$

\

We have $\dist(a, \ov R(\z_0)) < r$. Because 
supp $\mu=\partial A$ and  $\e$ and $r$ can be taken 
arbitrarily close to 0, this proves the density of  periodic 
points in  $\partial A$.

\

\

{\bf Section 2. Geometric coding trees , the complement of the 
proof of Theorem A.}

\

We shall prove a more abstract and general 
version of Theorem A 
here.  This will allow immediately to deduce Theorem A in the parabolic and non simply connected cases.

\

Let $U$ be an open connected subset of the Riemann sphere $\ov\C$.
Consider any holomorphic mapping $f:U\to \ov\C$ such that $f(U)\supset U$ and $f:U\to f(U)$ is a proper map. 
Denote $\Crit(f)=\{z:f'(z)=0\}$. This is called the set of 
critical  points for $f$. Suppose that $\Crit(f)$ is finite.
 Consider any $z\in f(U)$. Let $z^1,z^2,...,z^d$ be some
of the $f$-preimages of $z$ in $U$ where $d\ge 2$. Consider smooth curves $\g^j:[0,1]\to f(U)$, \ $j=1,...,d$,   joining $z$ to $z^j$ respectively (i.e. $\g^j(0)=z, \g^j(1)=z^j$), such that there are no critical values for iterations of $f$ in $\bigcup _{j=1}^d \g^j $, i.e. 
$\g^j\cap f^n(\Crit(f))=\emptyset$ for every $j$ and $n>0$. 

Let $\S^d:=\{1,...,d\}^{\Z^+}$ denote the one-sided shift space and $\s$ the shift to the left, i.e. $\s((\a_n))=(\a_{n+1})$. For every sequence $\a=(\a_n)_{n=0}^\infty \in \S^d$ we define $\g_0(\a):=\g^{\a_0}$. Suppose that for some $n\ge 0$, for every $ 0\le m\le n$, and all $\a\in\S^d$, the curves $\g_m(\a)$ are already defined.
Suppose that for $1\le m\le n$ we have $f\circ \g_m(\a)=\g_{m-1}(\s(\a))$, and $\g_m(\a)(0)=\g_{m-1}(\a)(1)$.

 Define the curves $\g_{n+1}(\a) $ so that the previous equalities hold by taking respective $f$-preimages of curves $\g_n$. For every  $\a\in\S^d$ and $n\ge 0$ denote $z_n(\a):=\g_n(\a)(1)$. 

For every $n\ge 0$ denote by $\S_n=\S^d_n$ the space of all 
sequences of elements of $\{1,...,d\}$ of length $n+1$. 
Let $\pi_n$ denote the projection $\pi_n:\S^d\to 
\S_n$ defined by $\pi_n(\a)=(\a_0,...,\a_n)$. As 
$z_n(\a)$ and $\g_n(\a)$ depends only on 
$(\a_0,...,\a_n)$, we can consider $z_n$ and $\g_n$ as functions on $\S_n$.

The graph $\T(z,\g^1,...,\g^d)$ with the vertices $z$ and $z_n(\a)$ and edges $\g_n(\a)$ is called a {\it geometric coding tree} with the root at $z$. For every $\a\in\S^d$ the subgraph composed of $z,z_n(\a)$ and $\g_n(\a)$ for all $n\ge 0$ is called a  {\it geometric branch} and denoted by $b(\a)$. The branch $b(\a)$ is called {\it convergent} if the sequence $\g_n(\a)$ is convergent to a point in $\cl U$. We define the {\it coding map} $z_\infty :{\cal D}(z_\infty)\to \cl U$ by $z_\infty(\a):=\lim_{n\to\infty}z_n(\a)$ on the domain ${\cal D}={\cal D}(z_\infty)$ of all such $\a$'s for which $b(\a)$ is convergent.

(This convergence is called in [PS] strong convergence.
In previous papers [P1], [P2], [PUZ] we considered mainly  convergence  in the sense $z_n(\a)$ is 
convergent to a point, but here we shall need the convergence of the edges $\g_n$.)

\

In the sequel we shall need also the following notation:  for each geometric branch $b(\a)$ denote by 
$b_m(\a)$ the part of $b(\a)$ starting from $z_m(\a)$ i.e.
consisting of the vertices $z_k(\a), k\ge m$ and of the 
edges $\g_k(\a), k>m$.

\

The basic theorem concerning convergence of 
geometric coding trees  is the following

\

{\bf  Convergence Theorem.} 1. Every branch except branches in a set of Hausdorff dimension  0 in a standard metric on $\S^d$, is 
convergent. (i.e $\HD(\S^d\setminus{\cal D} )=0$).
In particular for every Gibbs measure $\nu_\f$ for a 
H\"older continuous function $\f:\S^d\to \R$ \ 
$\nu_\f(\S^d\setminus{\cal D} )=0$, so the measure 
$(z_\infty)_*(\nu_\f)$ makes sense.

2. For every $z\in\cl U$ \ $\HD(z_\infty^{-1}(\{z\}))=0$. 
Hence for every $\nu_\f$ we have for the entropies: 
$\h_{\nu_\f}(\s)=\h_{ (z_\infty)_*(\nu_\f)}(\ov f)> 0$,
(if we assume that there exists $\ov f$ a 
continuous extension of $f$ to $\cl U$).

\

The proof of this Theorem can be found in 
[P1] and [P2] under some stronger assumptions (a slow convergence of $f^n(\Crit(f)$ to $\g^i$ for $n\to\infty$)  
To obtain the above version one should rely also on [PS] (where even $f^n(\Crit(f))\cap\g^i\not=\emptyset$ is allowed).

\

Recently, see [P4], a complementary fact was proved for 
$f$ a  rational map on the Riemann sphere, $U$ a 
completely invariant basin of attraction to a sink or a 
parabolic periodic point, under the condition (i) (see statement of Theorem C):

\

3.  Every $f$-invariant probability ergodic measure 
$\mu$, of positive  Lyapunov exponent, supported by 
$\cl z_\infty({\cal D})$ is a $(z_\infty)_*$-image of a 
probability $\s$-invariant measure on $\S^d$, (provided 
$f$ extends holomorphically to a neighbourhood of 
$\supp \mu$).

\

Suppose in Theorems B, C which follow,  that the  map $f$ extends 
holomorphically to a neighbourhood of  the closure of the limit set $\La$ of a tree , $\La =z_\infty({\cal D}(z_\infty))$. Then $\La$ is called a {\it quasi-repeller}, see [PUZ]. 

\

{\bf Theorem B.} For every quasi-repeller $\La$ for  
a geometric coding tree 

\noindent$\T(z,\g^1,...,\g^d)$ for a 
holomorphic  map $f:U\to\ov\C$, for every Gibbs measure $\nu$ for a H\"older continuous function $\f$ on $\S^d$ periodic points in $\La$ for the extension of $f$ 
to $\La$ are dense in $\supp (z_\infty)_*(\nu)$.

\

This is all we can prove in the general case. In the next Theorem  we shall introduce additional assumptions.

Denote 
$$\hat\La:=\{\hbox{all limit points of the sequences} \ z_n(\a^n), \a^n\in\S^d, n\to\infty\}$$

\

{\bf Theorem C.} Suppose we have a  tree as in 
Theorem B which 
satisfies  additionally the following conditions for every $j=1,...,d$: 
$$
\g^j\cap \cl(\bigcup_{n\ge 0} f^n(\Crit f)) = \emptyset ,
 \eqno {\rm {(i)}}
$$

There exists a neighbourhood $U^j\subset f(U)$ of $\g^j$ 
such that
$$
\Vol(f^{-n}(U^j) \to 0    \eqno {\rm {(ii)}}
$$
where Vol denotes the standard Riemann measure on $\ov\C$.

Then periodic points in $\La$ for $\ov f$ 
are dense in $\hat\La$.

\

\

Theorem C immediately follows from Theorem B if we 
prove the following: 

\

{\bf Lemma 3.} Under the assumptions of Theorem C 
(except we do not need to assume $f$ extends to $\ov f$)
for every  Gibbs measure $\nu$ on $\S^d$ we have  
$\supp (z_\infty)_*(\nu) = \hat\La$.

\

{\bf Proof of Lemma 3.} The proof is a minor modification of the proof of Convergence Theorem, part 1, but for the 
completness we give it here.

 Let $U^j$ and $U'^j$ be open 
connected simply connected neighbourhoods of $\g^j$ 
for $j=1,...,d$ respectively, such that $\cl U'^j\subset 
U^j$, \ $U^j \cap \cl(\bigcup_{n> 0} f^n(\Crit f)) = \emptyset $ and (ii) holds.

By (ii)  $\e(n):=\Vol(f^{-n}(\bigcup_{j=1}^d 
U^j)) \to 0$ as $n\to\infty$.

Define $\e'(n)=\sup_{k\ge n} \e(n)$. We have 
$\e'(n)\to 0$.

Denote the components of $f^{-n}(U^j)$ and of $f^{-
n}(U'^j)$ containing $\g_n(\a)$ where $\a_n=j$, by $U_n(\a), U_n'(\a)$  
respectively . Similarly to $z_n(\a)$ and $\g_n(\a)$ 
  each such component depends 
only on the first $n+1$ numbers in $\a$ so in our notation 
we can  replace $\a$ by 
$\pi_n(\a)=(\a_0,...,\a_n)\in\S_n$.

Fix arbitrary $n\ge 0$, \ $\a\in\S_n$ and $\d>0$. For 
every $m>n$ denote 
$$
\B(\a,m)=\{(j_0,...,j_m)\in\S_m: \ j_k=\a_k \hbox{\ \ for \ \ } k=0,...,n\}
$$
and 
$$\B_\d(\a,m)=\{(j_0,...,j_m)\in\B(\a,m):
 \Vol(U_m(j_0,...j_m)) \le \e(m)\exp(-(m-n)\d)\}.
$$
Denote also $\B(\a)=\pi_n^{-1}(\{\a\})  \subset \S^d$

\

Because all $U_m(j_0,...,j_m)$ are pairwise disjoint
$$
\sharp\S_m-\sharp \B_\d(\a,m) \le \exp (m-n)\d. 
 \eqno (7)
$$

By Koebe distortion theorem for the branches $f^{-m}$ 
leading  from $U^j\to U_m(\b)$ for $\b\in \S^d, \b_m=j$ 
we have 
$$
\diam(\g_m(\b)) \le \diam(U'_m(\b)) \le \Const(\Vol 
(U'_m(\b)))^{1/2}\le \Const(\Vol (U_m(\b)))^{1/2}
$$

Thus if $\b\in \B(\a)$ and $\pi_m(\b)\in \B_\d(\a,m)$ 
for every $m>m_0\ge n$ then for the length $b_{m_0}$ 
we have 
$$
\length (b_{m_0}(\b)) \le \Const \sum _{m>m_0}\e(m)^{1/2} 
\exp -(m-n)\d/2
$$

Now we shall rely on the following property of the 
measure $\nu$ true for the Gibbs measure for every 
H\"older continuous function $\f$ on $\S^d$:

There exists $\theta >0$ depending only on $\f$ such 
that for every pair of integers $k<m$ and every 
$\b\in\S^d$ 

$${\nu(\pi_m^{-1}(\pi_m(\b))) \over \nu(\pi_k^{-1}(\pi_k(\b)))}  < \exp -(m-k)\theta
$$

So  with the use of (7) we obtain
$$
{\nu(\B(\a) \setminus \bigcap_{m > m_0}\B_\d(\a,m)) 
\over
\nu(B(\a))}
 \le \sum_{m > m_0}   \exp (m-n)\d\exp(-(m-
n)\theta).
$$
We consider  $\d< \theta$.

\

As the conclusion we obtain the following 

{\bf Claim.}  For every $r>0, 0<\lambda<1$ if 
$n$ is large enough then for every $\a\in \S^d_n$ there 
is $\B'\subset\B(\a)$ such that
$$
{\nu(\B') \over \nu (\B(\a))} >\lambda
$$
and for every $\b\in \B' $
$$\length (b_n(\b)) < r.
$$

Indeed, it is sufficient to take $\B'=\bigcap_{m>m_0} 
\B_\d(\a,m)$ , where $m_0$ is the smallest integer $\ge 
n$ such that $\sum_{m>m_0}\exp(m-n)(\d-\theta)\le 1-\lambda$. (Of course the constant $m_0-n$ does not depend on $n,\a$.) Then for every $\b\in \B'$
$$
\length (b_n(\b)) < (m_0-n)\e'(n) + \Const(\e'(m_0))^{1/2}
\sum_{m>m_0}\exp(-(m-n)\d/2)   <  r
$$
if $n$ is large enough.

\

The above claim immediately proves our Lemma 3. $\hfill\clubsuit$

\

{\bf Remark   4.} \  Lemma 3 proves in particular (under the assumptions (i) and (ii) but without assuming $f$ extends to $\ov f$ ) that
$\cl\La =\hat\La$.

\

{\bf Remark 5.} \ Observe that Lemma 3 without any additional assumptions about the tree, like (i), (ii), is false.  For example take 
$z=p$ our sink, $z^1=p, z^j\not=p$ for $j=2,3,...,d$ and 
$\g^1\equiv p$. 
Then $p\in \La$ but $p\notin \supp (z_\infty)_*(\nu)$ 
for every Gibbs $\nu$. 

\

Observe  that if  
and (i) and (ii) are skipped in the assumptions of 
Theorem C then its assertion on the density 
of $\Lambda$ or the density of periodic 
points in $\hat\La$ is also false. We can take $z$ in a Siegel 
disc $S$ 
but $z$ different from the periodic point in $S$, $z^1\in 
S, z^j\notin S$ for $j=2,...,d$.

Here $\La$ is not dense even in the set $\La'$  intermediary between $\La$ and $\hat\La$ 
$$
\La':=\bigcup_{\a\in\S^d}\La(\a) \ \ \ \hbox{where}\ \ \ \La(\a):=
\{\hbox{the set of limit points of } z_n(\a), n\to\infty \}$$
(because $\La'$ contains a "circle" in the Siegel disc). 

$\hat\La$ corresponds to the union of impressions of all
 prime ends and $\La'$ corresponds to the union of all sets of principal points. See [P3] for this analogy.

\

We do not know whether Lemma 3 or Theorem C are  true without the assumption (i), only with the assumption (ii).

\

Now we shall prove Theorem B:

\

{\bf Proof of Theorem B.} We repeat the same scheme   as in Proof of Theorem A, the case 
discussed in Section 1.  Now $(\partial\D, g, l)$ is replaced 
by $(\S^d, \s, \nu)$. Its natural extension is denoted by $(\tilde\S^d, 
\tilde\s , \tilde\nu)$ \  (in fact $\tilde\S^d =\{1,...,d\}^{\Z}$). 
As in Section 1 we find a set $\tilde K$ with 
$\tilde\nu(\tilde K) > 1-2\e$ so that all points of 
$\tilde K$ satisfy the assumptions of Lemma 1 with constant $C, r$. The map $\ov R$ is 
replaced by $z_\infty$ and $Y$ is $\cl\La$ now.

Condition (1) makes sense along branches (which play the role of cones) , i.e. it takes 
the form:

there exists $M=M(r)$ arbitrarily large such that for every 
$\a\in\tilde K$

$$
b_M(\a)  \subset     B( z_\infty(\a), r/3).  \eqno (8)
$$

The crucial property we need to refer to Lemma 1 is 
$\chi_{(z_\infty)_*(\nu)}(\ov f) > 0$.  It holds because 
by Convergence Theorem, part 2, we know  that 
$\h_\nu(\s)=\h_{(z_\infty)_*(\nu)}(\ov f) > 0$ and by 
[R] 
$\chi_{(z_\infty)_*(\nu)}(\ov f) \ge  {1 \over 2}  \h_{(z_\infty)_*(\nu)}(\ov f)  > 0$

As in Section 1. for every  $\a=(...\a_{-1},\a_0,\a_1,...)  \in \tilde K$ there exists $N$ arbirarily large such that 
$\b=\pi_0\tilde\s^{-N}(\a)\in \tilde K$ is close to $\a$. In
 particular 
$$\b=(\a_0,\a_1,...,\a_M,w,\a_0,\a_1,...)
$$
where $w$ stands for a sequence of $N-M-1$ symbols from 
$\{1,...d\}$ and $N >  M$.

By (8) we have 
$$b_M(\a)  \subset  B(z_\infty(\a), r/3) \ \ \ \ \ \  \and$$
$$b_M(\b)   \subset  B(z_\infty(\b), r/3)$$
We have also 
$$z_M(\a)=z_M(\b).$$

So $\g:=\bigcup_{n=M+1}^{N+M}   \g_n(\b) \subset B(z_\infty(\a) ,r)$.  Since $F_N(z_\infty(\a))=z_\infty(\b)$ we 
have similarly as in Section 1, (6), 
$F_N(z_M(\a)=z_{M+N}(\b)$, i.e. $F_N$ maps one end of $\g$ to the other.  We have also, similarly to (5), 
$F_N(B(z_\infty(\a)), r)\subset F_N(B(z_\infty(\a)), r)$ and $F_N$ is a contraction. 

 One end of the curve $\Gamma$ built from $\g, F_N(\g), F^2_N(\g),...$
is periodic of period $N$, is in $B(z_\infty(\a),r)$ and 
is the limit of the branch of the periodic point 
$$
(\a_0,...,\a_M,w, 
\a_0,...,\a_M,w, ...) \in \S^d.
$$
  Theorem B is proved.\hfill $\clubsuit$

\

{\bf Proof of Theorem A. The conclusion.}

Denote $\Crit^+:=\bigcup_{n>0}f^n(\Crit(f)|_A)$.
Let $p$ denote the sink in $A$ or a parabolic 
point in $\partial A$. 

Take an arbitrary point $z\in A\setminus \Crit^+,\ z\not= 
p$. Take an arbitrary geometric coding tree 
${\cal T} (z,\g^1,...,\g^d) $ in $A\setminus (\Crit^+\cup 
\{p\})$, where  $d=\deg f|_A$.

Observe that (i) is satisfied because 
$\cl\Crit^+ =\{p\}\cup \Crit^+$.

Condition (ii) also holds because taking $U^j\subset A$ 
we obtain $f^{-n}(U^j)\to \partial A$, hence there exists 
$N>0$ such that  for every $n\ge N$ we have 
$$
f^{-n}(U^j)\cap U^j= \emptyset. 
$$ 

Indeed if we had $\Vol f^{-n_t}(U^j) > \e >0$ for a 
sequence $n_t\to\infty$ we could assume that $n_{t+1}-
n_t\ge N$. We would have $f^{-n_t}(U^j) \cap f^{-n_s}(U^j)=\emptyset$ 
for every $t\not= s$ hence

\noindent $\Vol\bigcup_t f^{-n_t}(U^j) = \sum_t\Vol f^{-n_t}(U^j) 
\ge \sum_t \e =\infty$, a contradiction. 

Thus we obtain from Theorem C that periodic points in 
$\La$ are dense in $\hat\La$.  The only thing to be 
checked is 

$$\hat\La = \partial A \eqno (9)$$

\

(If $A$ is completely invariant then $Z=\bigcup_{n\ge 0} 
f^{-n}(z)$ is a subset of $A$. It is dense in Julia set, in 
particular in $\partial A$. However in general situation $Z\not\subset A$ so the existence of a sequence in $Z$ 
converging to a point in $\partial A$ does not imply 
automatically the existence of such a sequence in $Z\cap A$.)

It is not hard to find a compact set $P\subset A$ such that  
$P\cap (\Crit^+\cup\{p\})= \emptyset $ and such 
that  for every $\z_0\in \partial A \setminus \{p\}$ for 
every $\z\in A$ close enough to $\z_0$, there exists 
$n>0$ such that $f^n(\z)\in P$. The closer $\z$ to 
$\z_0$, the larger $n$.

Cover $P$ by a finite number of topological discs $D_\tau \subset A$. 
There  exist topological discs $D'_\tau$ which union also covers $P$ such that $\cl D'_\tau \subset D_\tau$. Join each disc $D_\tau$ with $z$ by a curve $\d_\tau$ without selfintersections
disjoint with $\Crit^+$ and $p$. Then for every $\tau$ 
there exists a topological disc $V_\tau\subset A$  being a 
 neighbourhood of $D_\tau\cup \d_\tau$ also
disjoint with $\Crit^+$ and $p$. 

For every $\e>0$ there exists $n_0>0$ such that for 
every $n>n_0$ and every branch $F_n$ of $(f|_A)^{-n}$ 
on $V_\tau$
$$
\diam (F_n(D'_\tau\cup\d_\tau) ) <\e$$
by the same reason by which $\Vol f^{-n}(U^j) \to 0$ and 
next (by Koebe distortion theorem, see Proof of Lemma 3) 
$\diam f^{-n}(U'^j) \to 0$.

So fix an arbitrary $\z_0\in \partial A\setminus \{p\}$ 
and take $\z\in A$ close to $\z_0$. Find $N$ and 
$\tau$ such that $f^N(\z)\in D'_\tau$. We can assume $N>n_0$. Let $F_N$ be the branch of $f^{-1}$ on $V_\tau$ such that $F_N(f^n(\z))= 
\z$.  Then $\dist (\z, F_N(z)) <\e$. But $F_N(z)$ is a 
vertex of our tree. Letting $\e\to 0$ we obtain (9)\hfill$\clubsuit$

\

{\bf Remark 6.} \  One can apply Theorem C to $f$ a 
rational mapping on the Riemann sphere and $d=\deg 
(f)$ under the 
assumptions that for the Julia set $J(f)$ we have 
$\Vol J(f)=0$ and that the set $\cl\Crit^+$ does not dissect 
$\ov\C$. Indeed in this case we take $z$ in the immediate basin of a sink or a parabolic point 
and curves $\g^j$ disjoint with $\cl\Crit^+$. Then the 
assumptions (i), (ii) are satisfied, so periodic points in 
$\La$ are dense in $\hat\La$. A basic property of $J(f)$ 
says that 
$\bigcup_{n>0}f^{-n}(z)$ is dense in $J(f)$, i.e. $\hat\La =Jf)$, hence periodic points in $\La$ are dense in $J(f)$.

In this case however we can immediately deduce 
the density of periodic sources belonging to $\La$ in 
$J(f)$ from the fact that periodic sources are dense in $J(f)$ and from the theorem saying that every periodic 
source $q$ is a limit of a branch $b(\a), \ \a\in\S^d$ 
 converging to it. So $q$ belongs to $\La$ automatically. 
For details see [P4].

 \

\

\

\

\

{\bf References:}

\

[D] A. Douady, Informal talk at the Durham Symposium, 1988.

\

[Du] P. L. Duren,  Theory of $H^p$ spaces, New York -- 
London, Academic Press 1970.

\

[EL] A. E. Eremenko,  G. M. Levin,  On periodic points of polynomials. Ukr. Mat. Journal 41.11 (1989), 1467-1471.

\

[K] A. Katok, Lyapunov exponents, entropy and periodic 
points for diffeomorphisms,  Publ. Math. IHES 51 (1980), 
137-173.

\

[LP] G. Levin, F. Przytycki, External rays to periodic 
points, preprint.

\

[Pe]  C. L. Petersen, On the Pommerenke-Levin-Yoccoz 
inequality, preprint IHES M/91/43.

\

[P1] F. Przytycki, Hausdorff dimension of harmonic 
measure on the boundary of an attractive basin for a 
holomorphic map.  Invent. Math. 80 (1985), 161-179.

\

[P2]  F. Przytycki, Riemann map and holomorphic 
dynamics. Invent. Math. 85 (1986), 439-455.

\

[P3] F. Przytycki, On invariant measures for iterations of 
holomorphic maps. In "Problems in Holomorphic 
Dynamics", preprint IMS 1992/7,  SUNY at Stony Brook.

\

[P4] F. Przytycki, Accessability of typical points for 
invariant measures of positive Lyapunov exponents for 
iterations of holomorphic maps.

\

[PUZ] F. Przytycki, M. Urbanski, A. Zdunik,  Harmonic, Gibbs and Hausdorff measures for holomorphic maps. Part 1 in Annals of Math. 130 (1989), 1-40. Part 2 in Studia Math. 97.3 (1991), 189-225. 

\

[PS] F. Przytycki, J. Skrzypczak, Convergence and pre-images of limit points for coding trees for iterations of holomorphic maps, Math. Annalen 290 (1991), 425-440.

\

[R] D. Ruelle, An inequality for the entropy of 
differentiable maps, Bol. Soc. Bras. Math. 9 (1978), 83-87.

\

[Ro] V. A. Rohlin, Lectures on the entropy theory of 
transformations with invariant measures, Usp Mat.Nauk 
22.5 (1967), 3-56, (in Russian);  \ Russ. Math.Surv. 22.5 
(1967), 1-52.

\end